\documentclass[11pt, a4paper, twoside]{article}

\usepackage[utf8]{inputenc}
\usepackage{amsmath, amssymb}
\usepackage{amsfonts}
\usepackage{amsthm}
\usepackage[margin = 3 cm]{geometry}
\usepackage[colorlinks=true, citecolor=black, linkcolor=black, urlcolor=black]{hyperref}
\usepackage{graphicx}
\usepackage[dvipsnames]{xcolor}
\usepackage{float}
\usepackage{verbatim}
\usepackage{multirow}
\usepackage{geometry}
\usepackage{subcaption}
\usepackage{enumitem}
\usepackage{booktabs}
\usepackage[round]{natbib}
\usepackage{mathtools}
\mathtoolsset{showonlyrefs}
\usepackage{fancyhdr}

\theoremstyle{definition}

\newtheorem{proposition}{Proposition}
\newtheorem{example}{Example}[section]

\DeclareMathOperator*{\argmax}{\arg\!\max}

\usepackage{setspace}

\title{Bayesian Selective Inference: Non-informative Priors}
\author{Daniel G. Rasines and G. Alastair Young \\
Department of Mathematics, Imperial College London
}

\date{}

\begin{document}

\maketitle

\begin{abstract}
We discuss Bayesian inference for parameters selected using the data. First, we provide a critical analysis of the existing positions in the literature regarding the correct Bayesian approach under selection. Second, we propose two types of non-informative priors for selection models. These priors may be employed to produce a posterior distribution in the absence of prior information as well as to provide well-calibrated frequentist inference for the selected parameter. We test the proposed priors empirically in several scenarios.

\end{abstract}

\section{Introduction}

Selective inference considers problems in which the parameter of interest is selected using the data. It is well known that failing to acknowledge this adaptivity in the subsequent inferences yields the reported error assessments invalid, if the same data that was used for selection is then used for inference. For example, type I error guarantees of frequentist testing procedures are lost. In this work we review selective inference from a Bayesian viewpoint and address the problem of specifying non-informative priors in these situations.

The settings considered here can be formalised as follows. Suppose we have data $Y$ whose distribution has a density or probability mass function $f(y; \theta), y\in \mathcal Y$, known up to a finite-dimensional parameter $\theta\in \Theta$, and that there exists a potential subparameter of interest, $\psi \equiv \psi(\theta)$, that may be selected for future study after observing the data, possibly by an artificially randomised procedure. More precisely, we assume that there exists a function $p\colon \mathcal{Y}\to [0, 1]$ determined prior to the data collection such that, having observed $Y = y$, inference on $\psi$ is performed with probability $p(y)$. The function $p(y)$, which fully characterises the selection mechanism, will be referred to as the selection function. 

A typical example of selection problem is the so-called selected mean problem, where we have $Y = (Y_1,\ldots, Y_m)^T$ for some $m > 1$, with $Y_i\sim N(\theta_i, 1)$ independently and $\theta = (\theta_1, \ldots, \theta_m)^T \in \mathbb{R}^m$. In this setting the parameter of interest is typically defined as the subset of $\theta_i$'s that produce the largest $k$ observations for some prespecified $k < m$, or those that produce observations above a given threshold $t \in \mathbb{R}$ determined by some formal rule, such as a multiplicity adjustment of a significance test. This model and generalisations of thereof have applications in many areas. For example, in biostatistics each $\theta_i$ might represent a treatment effect for a given disease, and it could be the case that only the most promising treatments are selected for further investigation due to economical or time constraints.

Recently, selection issues have attracted much attention in the regression literature. Let $Y \in \mathbb{R}^n$ be a response vector and $X\in \mathbb{R}^{n\times p}$ be a known, fixed design matrix containing the observed values of $p$ covariates. Suppose that a variable-selection algorithm is used to select a non-empty subset $s\subseteq \{1, \ldots, p\}$ of the covariates, and that we wish to provide inference for $\psi = \beta$ in the model $Y\sim N(X(s)\beta, \sigma^2 I_n)$, where $X(s)$ is the submatrix of $X$ containing the observations of the selected covariates, $\theta = (\beta^T, \sigma^2)^T$, and $I_n$ is the $n\times n$ identity matrix. Alternatively, we might consider the less restrictive model $Y\sim N(\mu, \sigma^2 I_n)$, where $\mu \in \mathbb{R}^n$ is unconstrained, and provide inference for the best linear predictor of $Y$ in the selected model, defined as $\psi = \{X(s)^T X(s)\}^{-1} X(s)^T \mu$; see \cite{berketal}. In both cases the selection function would be $p(y) = \mathbf{1}\{y\in E(s)\}$, where $E(s)\subseteq \mathbb{R}^n$ is the set of observations for which the variable-selection algorithm would have produced the same set of selected covariates, and $\mathbf{1}\{\cdot\}$ is the indicator function. 

From a repeated-sampling viewpoint, selection alters the sampling distribution of the data by favouring data points with a higher selection probability. This has clear implications in the frequentist paradigm and has led to the development of the so-called conditional approach, which advocates basing inference for the selected parameter on the conditional distribution of the data given selection. This approach was lucidly formalised by \cite{fithianetal} and framed withing Fisherian statistical thinking by \cite{kuffneryoung}. The Bayesian standpoint regarding selection is less clear. The classical view is that inference should not be altered by selection. The argument is that, since Bayesian inference operates conditionally on the data, in particular it conditions on the selection event, which is therefore automatically accounted for, as explained by \cite{dawid}. This viewpoint was however questioned by several authors, such as \cite{yekutieli}, who argued that in some situations the posterior distribution has to be appropriately modified in the presence of selection.

In the following section we analyse this discrepancy between the frequentist and Bayesian approaches to inference under selection. We argue that, in general, an adjustment for selection is necessary in order to achieve approximate repeated-sampling validity. Formally, this adjustment is achieved by attaching the prior density to the likelihood of the conditional distribution of the data given selection. Then, in Section \ref{SEC: prior}, we consider the problem of specifying non-informative prior densities for selected parameters of exponential families via the requirement of accurate repeated-sampling calibration. These priors admit a simple expression and enable us to extend inferential methods for selective normal models to more general settings in a relatively simple manner. We propose the use of either a data-dependent prior, which by construction achieves the appropriate frequentist probability matching, or a Jeffreys prior, constructed from the selective likelihood, which is demonstrated to provide similar results and is easier to implement. The example of `play the winner', which involves a high-dimensional parameter, is considered in Section \ref{SEC: play winner}. This practically important example has been a key element of discussion on Bayesian inference under selection; see \cite{dawid} and \cite{senn}.

\section{Bayes and selection} \label{SEC: bayes}

Different approaches to statistical inference lead to two opposing views as to the correct analysis of the data in the presence of selection. On the one hand, frequentist methods evaluate the accuracy of inferential procedures with respect to the sampling distribution of the data at a fixed value of the parameter. Since selection modifies the sampling distribution by favouring data points with higher selection probability, it is clear that the reported accuracy should be appropriately modified. On the other hand, Bayesians typically adopt the view that, once the data has been observed, the recognition that a different realisation could have resulted in a different inferential problem, or in no problem at all, should have no effect on the inference \citep{dawid}. 

Adoption of the first viewpoint leads to adherence to the so-called conditional approach, which advocates that inference for the selected $\psi$ should be based on the conditional distribution of the data given selection. Such distribution has density or mass function
\begin{equation} \label{EQ: conditionalmodel}
f(y\mid S; \theta ) = \frac{f(y; \theta) p(y)}{\varphi(\theta)}, \quad \varphi(\theta) = \mathbb{E}_\theta[p(Y)] ;
\end{equation}
where the normalising constant $\varphi(\theta)$ is the probability that $\psi$ gets selected when $\theta$ is the true parameter, and $S$ denotes the selection event. The motivation for basing inference on \eqref{EQ: conditionalmodel} is that, under repeated sampling from a given $f(y; \theta)$, if inferences are only provided for those samples that get selected, the reported error assessments are well-calibrated. For example, nominal $1-\alpha$ confidence sets would contain the true parameter at least $(1-\alpha)\%$ of the times they are reported. Henceforth we will refer to the conditional distribution of $Y$ given selection as the selective distribution, and to the corresponding likelihood function, $L_S(\theta) = f(y\mid S; \theta )$, as the selective likelihood.

An intuitive way to interpret the conditional approach is as a form of information splitting. For a given $\psi$, let $R$ be the random variable that takes the value $1$ if $\psi$ gets selected and $0$ otherwise. That is, $R\vert Y \sim \mathrm{Bernoulli}\{p(Y)\}$. Following \cite{fithianetal}, the data-generating process of $Y$ may be thought of as consisting of two stages. In the first stage, the value of $R$, $r$ say, is sampled from its marginal distribution, and in the second, $Y$ is sampled from the conditional distribution $Y \mid r$. Since it is $R$ that determines whether we are going to provide inference for $\psi$ or not, inference based on information revealed in stage two is necessarily free of any selection bias, as it removes the information about the parameter provided by $R$. 

In the Bayesian literature the appropriate mode of inference under selection is less clear. The predominant viewpoint until recent years was that Bayesian inference should not be modified in the presence of selection. Quoting \cite{dawid}:
\begin{quote} 
\textit{Since Bayesian posterior distributions are already fully conditioned on the data, the posterior distribution of any quantity is the same, whether it was chosen in advance or selected in the light of the data: that is, for a Bayesian, the face-value approach is fully valid, and no further adjustment for selection is required.}
\end{quote}
Expressed differently, the fact that the statistician decides to focus on certain aspects of the data after it has been observed does not alter its sampling distribution, which should therefore remain unchanged. 

As Dawid points out, the contrast between the Bayesian and the frequentist standpoints is somewhat paradoxical. In many situations, Bayesian analyses formally match, either exactly or approximately, face-value frequentist analyses, and it is universally agreed that the latter methods are not valid in the presence of selection. Why, then, would such results be correct if reached by a Bayesian argument? Dawid argues that, certainly, whenever a frequentist approach is unreasonable, so is any Bayesian approach that provides similar answers. However, according to Dawid, this does not reflect a fundamental issue about the Bayesian updating mechanism under selection. Rather, it is the result of a poor prior specification.

To understand why some prior distributions can be problematic under selection, consider the example, due to \cite{dawid}, of providing inference for the mean of $Y\sim N(\theta, 1)$ only if $Y > 0$, and consider the standard class of conjugate priors for $\theta$ given by $\{N(0, \lambda^2)\colon \lambda > 0\}$. For a given $\lambda$, the posterior distribution of $\theta$ given $Y = y$ is normal with mean $\{\lambda^2 / (1 + \lambda^2)\} y$ and variance $\lambda^2 / (1 + \lambda^2)$. Small values of $\lambda$ produce a shrinkage effect, pulling the posterior density of $\theta$ towards $0$, while large values produce posterior inferences which are very similar to those provided by a face-value frequentist approach. Now, consider this problem from a frequentist perspective. If the true $\theta$ is significantly greater than zero, the selection bias will be very small and almost no adjustment is needed. Conversely, if $\theta$ is small, selection needs to be accounted for, as $Y$ will tend to overestimate it. Therefore, the Bayesian analysis described before carried out with a large value of $\lambda$ is appropriate only if we expect $\theta$ to be large, and will lead us astray if it is not. This is obvious, in a way: we would not adjust for selection if we believed a priori that $\theta$ is large and that an adjustment would not be necessary, and we would adjust for it if we suspected that $\theta$ is small and that therefore an adjustment would be required.   

Of course, the fact that the prior choice is important and should not be taken lightly is not surprising. The key issue that the previous analysis highlights is that the impact of the prior on the analysis is stronger under selection, and, perhaps more importantly, that non-informative or weakly informative Bayesian analyses in these contexts are much more challenging than in non-selective settings. In the previous problem, if no selection takes place, lack of information about $\theta$ can be dealt with by taking a large, or even infinite, value of $\lambda$, as this minimises the influence of the prior on the analysis. In the presence of selection, however, this is no longer appropriate, as such theoretically ``non-informative'' priors do in fact entail critical information about the parameter, namely, its likelihood to lie in a region of the parameter space for which a selection adjustment would be appropriate.

\begin{figure}[ht]
\centering
\includegraphics[scale = 1]{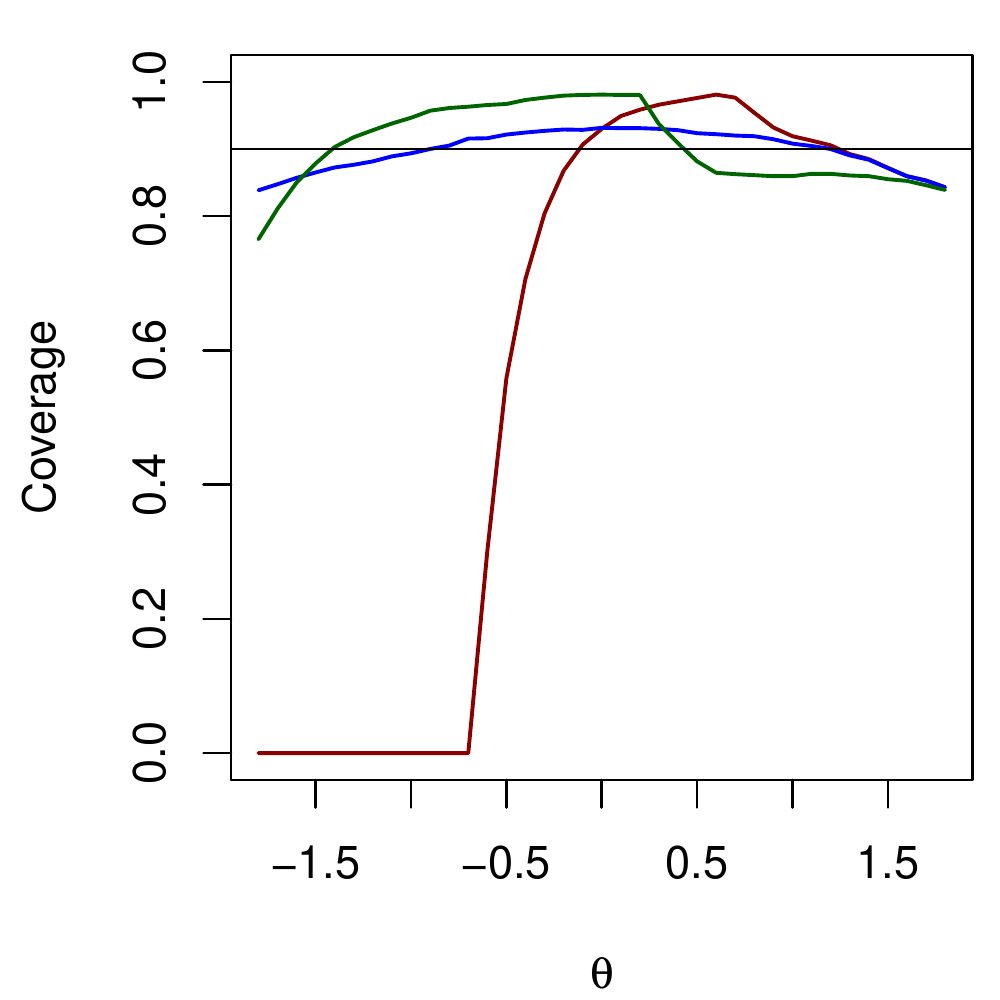}
\caption{Coverage of credible intervals under repeated sampling from fixed values of $\theta$ derived from the unadjusted posterior (the red line corresponds to sampling with selection and the blue one without it), and from the selective posterior (green).}
\label{FIG: coverages}
\end{figure}

However, even if we are comfortable with the implications of a certain prior in a selective problem, unadjusted Bayesian inference can still be problematic from a repeated-sampling viewpoint. Consider the model $Y\mid \theta \sim N(\theta, 0.2)$ and suppose that the prior distribution of $\theta$ is standard normal. In Figure \ref{FIG: coverages} we have plotted the observed coverages of equal-tailed $0.9$-credible intervals for fixed values of the true parameter in two different sampling regimes: a selective one, where only the observations with $Y > 0$ where kept, and a non-selective one, where we constructed the intervals for all the sampled values of $Y$. The curves where obtained by simulation. We see that in the selective regime the intervals are very poorly calibrated for small values of $\theta$, even for values with non-negligible selection probability (for example, when $\theta = -0.5$, the probability that $Y > 0$ is 0.13). By contrast, the coverage figures in the non-selective settings are much more stable. The instability observed in the selective regime is clearly undesirable, and, while a purist may argue that inferences are still well-calibrated once the coverage curve is averaged over the prior distribution, which in this case is $\pi(\theta\mid S)\propto \phi(\theta) \varphi(\theta)$, where $\phi$ is the standard normal density, most Bayesians would presumably agree that this behaviour is too extreme.  

To avoid issues relating to selection, a natural option would be construct the posterior distribution using only the information about $\theta$ provided by $Y\mid \{R = 1\}$, where as before $R$ is the indicator variable of the selection event. This way, the prior distribution is updated using only information which is not influenced by selection. Denoting the prior density by $\pi(\theta)$, the resulting posterior density would be
\begin{equation} \label{EQ: selective posterior}
\pi(\theta\mid y) \propto \pi(\theta) L_S(\theta) \propto  \frac{\pi(\theta) f(y;\theta) }{\varphi(\theta)} .
\end{equation}
In the remainder of this work we will refer to posteriors of this form as selective posteriors. Inference based on selective posteriors allows the injection of prior information while avoiding potential problems arising from selection.

These type of posteriors have been considered by several authors, typically in settings where there is a explicit bias in the sampling mechanism of the data, like survival models. They are discussed by \cite{bayarri-degroot} and \cite{bayarri-berger-selection}, and they also appear, in a different context, in \cite{bayarri-postpred}, where they are referred to as ``partial posterior densities" and are used for computing Bayesian $p$-values in a two-step procedure. As noted by \cite{harville}, posterior \eqref{EQ: selective posterior} can also be thought of as following from a standard Bayesian updating of the modified prior $\pi^*(\theta) = \pi(\theta)/\varphi(\theta)$. This could be given the following interpretation: since unadjusted inference is very imprecise for values of $\theta$ with small selection probability, we assign a larger prior probability to those values to achieve a higher protection in those regions. In Figure \ref{FIG: coverages} we have plotted the coverage curve produced by the selective posterior in the selective regime, which shows a much smoother behaviour than the unadjusted posterior.   

\subsection{Fixed and random parameters}
The belief that Bayesian inference is unaffected by selection has also been questioned by other authors on more conceptual grounds. \cite{MR1, MR2} considered a sequence of binomial experiments in which only some of the inferences, corresponding to promising results in the setting considered, were reported. They made the observation that, if the sequence of parameters analysed had been sampled independently from a prior distribution, then no correction for selection was needed, in agreement with the classical stance. However, if the sequence of parameters had been generated by sampling one parameter from the prior and then successively testing it, then the appropriate inference needed a correction for selection. This idea was later extended and refined by \cite{yekutieli}, who considered more general settings and linked these ideas to a Bayesian analogue of the false discovery rate.
 
According to Yekutieli, the correct Bayesian inference for a selected parameter depends on how the selection mechanism acts on the parameter space. Consider the joint sampling model of $(\theta, Y)$ and a selection function $p(y)$. Yekutieli calls $\theta$ \textit{random} if the joint sampling scheme for the parameter and data is such that pairs $(\theta, Y)$ are sampled from their joint distribution until $\theta$ gets selected, and calls $\theta$ \textit{fixed} if it is sampled from its marginal distribution, held fixed, and $Y$ is sampled from the conditional distribution $Y\mid \theta$ until selection takes place. 

As before, let $R$ be the binary random variable that indicates if selection of the parameter under consideration has happened. If $\theta$ is random, its density given selection is $\pi(\theta\mid R = 1) \propto \pi(\theta) \varphi(\theta)$, whereas if it is fixed, its conditional density in unchanged, $\pi(\theta\mid R = 1) = \pi(\theta)$. In the case of the data, its conditional density given $\theta$ and the selection event is $f_S(y; \theta) =  f(y; \theta) p(y)/\varphi(\theta)$ in both cases. Thus, the posterior distribution for a random parameter is 
\begin{equation} \label{EQ: random}
\pi(\theta\mid y) \propto \frac{\pi(\theta) \varphi(\theta) f(y;\theta) }{\varphi(\theta)} = \pi(\theta) f(y;\theta),
\end{equation}
while for a fixed parameter it is given by
\begin{equation} \label{EQ: fixed}
\pi(\theta\mid y) \propto  \frac{\pi(\theta) f(y;\theta) }{\varphi(\theta)} .
\end{equation}
Hence, Bayesian inference for a selected random parameter is unaffected by selection, but for a fixed parameter the posterior needs to be adjusted, and it is formally constructed by attaching the prior density to the selective likelihood. Yekutieli argues that Dawid's analysis follows from the implicit assumption that the parameter is random.

To elucidate the difference between these two regimes let us consider the following example, which is a slight variation of the first example in \cite{yekutieli}. Suppose that $\theta$ represents the average academic ability of the students in a given population. Each student takes an exam and gets a grade $Y_i$, and the students whose grade exceeds a given threshold $t$ are admitted to university. If a person in the university, who has access to the grades obtained by the admitted students and to prior knowledge about $\theta$ in the form of historical data, say, wants to estimate the average academic ability in the population, they should carry out the analysis treating $\theta$ as a fixed parameter. If, instead, $\theta$ is a student-specific parameter that represents the student's academic ability, and the person in the university wants to estimate each of the academic abilities of the admitted students, then each estimated $\theta$ should be treated as a random parameter. 

The previous analysis, however, does not cover situations with non-informative priors. It is argued by Yekutieli (Section 2.2) that, in cases where a non-informative prior would have been used without selection, a non-informative prior should also be used with selection, regardless of whether $\theta$ is random or fixed. The reason given is that, just as in the non-selective regime, all the information about the parameter in the selective problem comes from the data. Yekutieli further suggests to use the same non-informative prior in the selective and non-selective models, because the decision of providing inference for a particular parameter should not alter the prior.

We agree with the first assertion but not with the second one, that is, we believe that the non-informative prior used in the selective model should be different from the one that would be used in the non-selective one. Non-informative priors are typically designed to achieve certain properties with respect to the sampling model of the data. In a selection model, both in the random and fixed parameter regimes, the sampling distribution of the data is the same (it has density or mass function $f(y\mid S; \theta)$), so the non-informative prior should also be the same in both cases. Furthermore, since the sampling distribution of the data depends on the selection mechanism, we argue that any appropriate choice of non-informative prior should also depend on selection. Consider, for example, a location model $f(y;\theta) = g(y-\theta)$. The standard non-informative prior for $\theta$ in absence of selection is $\pi(\theta) \propto 1$, which enjoys several interesting properties: it is Jeffreys, minimises the prior influence, and produces inferences with a valid frequentist interpretation. However, all these nice properties no longer hold if inference for $\theta$ is only provided for certain values of $y$, as the selective model $\{g(y - \theta) p(y)/\varphi(\theta)\colon \theta\in \Theta\}$ is not a location model unless $p(y)\propto 1$. We expand on these ideas in the following section.

As for informative priors, while the posterior densities \eqref{EQ: random} and \eqref{EQ: fixed} are formally correct given the respective sampling mechanisms, it is not clear that a parameter can be labelled as random or fixed without explicit consideration of this sampling process, which may not be well-defined in some cases. We therefore suggest to use the adjusted posterior to avoid any potential bias that may arise due to selection. We note, however, that the distinction between random and fixed parameters is clarifying, in that it captures the two types of bias that can arise because of selection. The first occurs in the sample space, when the same parameter is analysed multiple times: if some of the analyses are not reported because the parameter did not appear to be significant in light of the data, the published reports about it will on average overstate this significance. The second one is a bias in the parameter space: parameters that are in some way interesting are more likely to be analysed than those which are not.

\section{Non-informative priors for selective inference} \label{SEC: prior}


Non-informative priors allow the derivation of posterior distributions without explicit incorporation of prior information. They serve different purposes: the resulting posterior can be employed as a reference against which posteriors derived from subjective priors can be compared; they can be used to derive frequentist methods that retain some of the appealing characteristics of Bayesian methods; and they allow to carry out an analysis when very little or no prior information is available. One natural way to derive non-informative priors is by requiring that they are approximately well-calibrated under repeated sampling from a fixed parameter. More precisely, that they satisfy
\begin{equation} \label{EQ: PMP condition}
\mathbb{P}_\theta\{ \psi \leq \Pi^{-1}(\alpha\vert Y)\} = \alpha + \varepsilon(\alpha, \theta),
\end{equation}
where $\Pi(\psi\vert Y)$ is the marginal posterior distribution function of the interest parameter and $\varepsilon(\alpha, \theta)$ is small for all $\alpha\in (0, 1)$ and $\theta\in \Theta$. Condition \eqref{EQ: PMP condition}  ensures that posterior claims about the parameter have approximate validity in a frequentist sense. For example, $(1-\alpha)$-credible intervals for $\psi$ are also approximate $(1-\alpha)$-confidence intervals. In the context of selective inference, this condition would be required to hold with respect to the selective distribution of the data.

In non-selective regimes involving independent and identically distributed observations from a regular model, the quantile-matching requirement holds with an error decreasing at rate $n^{-1/2}$ for any prior that places positive probability around the true parameter, where $n$ is the sample size. In addition, certain priors, known as probability-matching priors, improve on this error rate, lowering it to $n^{-1}$ in general, and to $n^{-3/2}$ for some distributions; see \cite{DattaMukerjee}. Unfortunately, the formal analysis leading to the derivation of these priors cannot be replicated in selection models, as it relies on the asymptotic normality on the posterior, which does not hold under selection except in trivial cases.

\cite{yekutieli} argues that non-informative priors should not be altered in the presence of selection. In particular, for a selective normal-location model, the improper uniform prior $\pi(\theta)\propto 1$ is suggested. This prior is also used by \cite{panigrahi-scalable} for the case where $\theta$ is a linear regression parameter. We argue that a non-informative prior \textit{should} be appropriately modified in the presence of selection, to account for the change of sampling distribution. As the following result illustrates, not doing so can lead to poorly calibrated posterior inferences. The proof can be found in Appendix \ref{APP: prop}.

\begin{proposition} \label{PROP: uniform prior}
Let $Y \sim N(\theta, \sigma^2)$, with $\sigma^2>0$ known, and $p(y) = \mathbf{1}(y>t)$ for some fixed $t\in \mathbb{R}$, and let
\begin{equation}
\Pi(\theta\vert y) = \frac{\int_{-\infty}^\theta \phi\{\sigma^{-1}(\tilde{\theta} - y) \} \varphi(\tilde\theta)^{-1} \mathrm{d}\tilde\theta}{\int_{-\infty}^\infty \phi\{\sigma^{-1}(\tilde{\theta} - y) \} \varphi(\tilde\theta)^{-1} \mathrm{d}\tilde\theta}
\end{equation}
be the selective posterior distribution function based on the uniform prior $\pi(\theta) \propto 1$. Then,
\begin{equation}
\mathbb{P}_{\theta_0}\{\theta_0 \leq \Pi^{-1}(\alpha\vert Y)\mid S\} < \alpha \quad \forall (\alpha, \theta_0) \in (0,1)\times \mathbb{R}.
\end{equation}
\end{proposition}

Proposition \ref{PROP: uniform prior} says that, when selection takes the form of one-sided truncation, the uniform prior produces inferences that overstate values of $\theta$ with low selection probability. So, in a way, correcting for selection overcompensates the selection bias. In addition, simulation results at the end of this section suggest that this claim holds for any selection mechanism with a increasing selection function $p(y)$.
In Section \ref{SEC: PMP_1d} we will see that this overcompensation also occurs in non-normal models under standard non-informative priors.
 
Unlike in the non-selective case, where the uniform prior leads to exact probability matching, it can be shown that in a selective normal-location model exact matching can only be achieved with a data-dependent prior. Before doing so, we introduce the concept of a confidence distribution. Consider the model $Y\sim N(\theta, \sigma^2)$, with $\sigma^2$ known, and an arbitrary selection function $p(y)$ which does not vanish almost everywhere. Conditionally on selection, the $p$-value function
\begin{equation}
H(\theta; y) = \mathbb{P}_\theta\left( Y \geq y \mid S \right)
\end{equation}
is uniformly distributed over $(0, 1)$ when $\theta$ is the true parameter value. Furthermore, if $y$ is such that $0 < H(0; y) < 1$, then $\theta\to H(\theta; y)$ is a distribution function. To see this, assume that $\sigma^2 = 1$ for simplicity and write 
\begin{equation}
H(\theta; y) = \frac{\mathbb{P}(S, Y\geq y)}{\varphi(\theta)} = \frac{\int_y^\infty \phi(\tilde y - \theta) p(\tilde y)\mathrm{d}\tilde y}{\int_{-\infty}^\infty \phi (\tilde  y - \theta) p(\tilde  y)\mathrm{d}\tilde  y}.
\end{equation}
Differentiating with respect to $\theta$ gives
\begin{equation}
\frac{\partial}{\partial \theta} H(\theta; y) = \frac{\mathbb{P}_\theta(S, Y\geq y)}{ \varphi(\theta)} \left\{ \mathbb{E}_\theta\left[ Y \mid S, Y\geq y \right] - \mathbb{E}_\theta\left[ Y\mid S\right] \right\} > 0
\end{equation}
for all $\theta$. Furthermore, $0 < H(0; y)$ implies that $\int_y^\infty \phi(\tilde y - \theta) p(\tilde y)\mathrm{d}\tilde y \geq K \{ \Phi(y_2 - \theta) - \Phi(y_1 - \theta) \}$ for some $y < y_1 < y_2 $ and $K > 0$, and we clearly have that $\mathbb{P}(S, Y< y) \leq \Phi(y - \theta)$. Therefore
\begin{equation}
\frac{1}{H(\theta; y)} = \frac{\varphi(\theta)}{\mathbb{P}(S, Y\geq y)} = 1 + \frac{\mathbb{P}(S, Y < y)}{\mathbb{P}(S, Y\geq y)} \leq 1 + \frac{1}{K} \frac{\Phi(y - \theta)}{\Phi(y_2 - \theta) - \Phi(y_1 - \theta) } \to 1
\end{equation}
as $\theta\to \infty$. An analogous argument, using that $H(0; y) < 1$, shows that the limit as $\theta \to -\infty$ is zero, so $H(\theta; y) $ is indeed a distribution function.

Considered as a random function of $Y$, $H(\theta; Y) $ is an example of a so-called confidence distribution. In a one-dimensional statistical model $\{f(y; \theta)\colon \theta\in \Theta\subseteq \mathbb{R}\}$, a confidence distribution is any function $G(\theta; Y)$ that satisfies the following two conditions:
\begin{enumerate}
\item $G(\theta; Y)$ is uniformly distributed over $(0, 1)$ when $\theta$ is the true parameter;
\item $\theta \to G(\theta; y)$ is a distribution function for every $y$.
\end{enumerate} 
Note that, in the case discussed here, $H(\theta; y)$ is a distribution function if $0 < H(\theta; y) < 1$ for some $\theta$, but since $H(\theta; Y)\sim U(0, 1)$ under $\theta$, it follows that $\mathbb{P}_\theta\{H(\theta; Y) \in (0, 1) \mid S\} = 1$, so the second condition is satisfied almost surely; note that the condition $0 < H(0; y) < 1$ trivially implies that $0 < H(\theta; y) < 1$ for all $\theta$.

Confidence distributions can be thought of as probability distributions over the parameter space that are induced by the data and can be derived without prior information. They are, in essence, a modern formulation of Fisher's fiducial probability \citep{fiducial}, and have attracted much attention in recent years due to their joint frequentist-Bayesian spirit; see \cite{xie} for a recent review. Confidence distributions can be used to derive confidence intervals with any desired level since, by definition, the interval $[H^{-1}(\alpha_1; Y), H^{-1}(\alpha_2; Y)]$ is a confidence interval with coverage $\alpha_2 - \alpha_1$ for any $0 < \alpha_1 < \alpha_2 < 1$. In addition, $H^{-1}(1/2; Y)$ is a median-unbiased estimator of $\theta$, as it satisfies $\mathbb{P}_\theta\{H^{-1}(1/2; Y) \leq \theta\mid S\} = 1/2$.

Given a confidence distribution $H(\theta; y)$, one can trivially construct an exact probability-matching prior by requiring that the posterior distribution function $\Pi(\theta\vert y)$ equals $H(\theta; y)$. This type of construction appears, for example, in \cite{reid_dataPMP}, and can be used as a basis for deriving probability-matching priors in more complex scenarios, as we will do in the next section. The resulting prior density is usually data-dependent, and is computed as 
\begin{equation}\label{EQ: PMP}
\pi_y(\theta) \propto -\frac{\frac{\partial}{\partial \theta} H(\theta; y)}{\frac{\partial}{\partial y} H(\theta; y)} .
\end{equation}
In the case of the selective normal-location model, this is, in addition, the only \textit{sensible} prior that provides exact matching. Indeed, consider the matching equation
\begin{equation}
\mathbb{P}_\theta\left \{\Pi(\theta\vert Y) \leq \alpha \vert S \right \} = \alpha, \quad \alpha\in (0, 1).
\end{equation}
Since the model is stochastically increasing in $\theta$, $\Pi(\theta\vert y)$ ought to be a decreasing function of $y$ for every $\theta$. Denoting its inverse with respect to $y$ by $l_\theta(\alpha)$, we have that
\begin{equation}
\mathbb{P}_\theta\left \{Y \geq l_\theta(\alpha) \vert S \right \} = H\{\theta; l_\theta(\alpha)\} = \alpha, \quad \alpha\in (0, 1).
\end{equation} 
Letting $\alpha = \Pi(\theta\vert y)$ gives the equality. 

For a general selection mechanism the functional form of $\pi_y(\theta)$ does not admit a simple expression. In what follows we will restrict the the practical considerations to a class of natural selection mechanisms following from one-sided truncation of the data. Let us assume that the data consists of two batches of random samples of sizes $n_1$ and $n_2$ from the distribution $N(\theta, 1)$, with respective sample means $Y_1\sim N(\theta, n_1^{-1})$ and $Y_2\sim N(\theta, n_2^{-1})$, and suppose that the selection criterion is $Y_1 > t$ for some threshold $t\in \mathbb{R}$. This situation arises when either the data has been artificially split in order to achieve higher inferential power, or when only the first batch is initially available and the second one is gathered after selection occurs in order to obtain more information about the parameter. We will in particular consider the case $n_2 = 0$, corresponding to a situation where all the data available for inference has been used for selection. 

Letting $n = n_1 + n_2$, $\gamma = n_1/n$, and $Y = \gamma Y_1 + (1 - \gamma) Y_2$, this model can be reduced by sufficiency to one involving a single observation $Y \sim N(\theta, n^{-1})$ and selection function $p(y) = \mathbb{P}( y + W > t)$, where $W \sim N(0, n^{-1}\{\gamma^{-1}-1 \})$ is independent of $Y$. To see this, simply define $W = (1-\gamma)( Y_1 - Y_2)$, which follows the claimed distribution and is independent of $Y$, as they are uncorrelated and normal, and note that $Y_1 = Y + W$. This analysis shows that, in a normal location model, selection based on a $\gamma$-split of the data is equivalent to an additive randomised selection scheme with randomisation noise $W$, as both selection strategies produce the same selection function.

Note that many `selected mean' problems can be reduced to a situation of the previous type via conditioning. Suppose that there are $m$ independent batches of size $n_1$ from $m$ different populations, with respective sample means $Y^{(i)}\sim N(\theta_i, n_1)$, and that inference is sought for the means yielding the largest $k$ observations, which we may assume to be $\theta_1, \ldots, \theta_k$ without loss of generality. If we condition on the data from the non-selected parameters, $Y^{(k+1)}, \ldots, Y^{(m)}$, the resulting selective model becomes
\begin{equation}
f(y^{(1)}, \ldots, y^{(k)} \mid S, y^{(k+1)}, \ldots, y^{(m)}; \theta) = \prod_{i = 1}^k f(y^{(i)}\mid Y^{(i)} > t; \theta_i)  
\end{equation}
where $t = \max\{y^{(k+1)}, \ldots, y^{(m)}\}$.
Thus, the problem can be written as $k$ separate inferential problems of the previous type. 

For these models, the exact probability matching prior can be written as
\begin{equation}
\pi_y(\theta) \propto 1 - \frac{ h_1\{n^{1/2}(\theta - t)\}}{h_1\{ n^{1/2}(\theta - y)\}} 
\end{equation}
when $n_2 = 0$ ($\gamma = 1$), where $h_1(x) = (\partial / \partial x) \log \Phi(x) = \phi(x)/\Phi(x)$. For $\gamma > 1$, a few calculations show that it can be simplified to
\begin{equation} \label{EQ: PMP_normal_long}
\begin{split}
& \pi_y(\theta) \propto \\ 
& \gamma^{1/2}\frac{\phi(n_1^{1/2}(\theta - t))}{\phi(n^{1/2}(\theta - y))} \\ 
 & \times \left\{ 1 - \Phi\left\{ \frac{n^{1/2}(y - \theta + \gamma(\theta - t)) }{(1 - \gamma)^{1/2}}\right\}-\frac{\int_y^\infty n^{1/2}\phi(n^{1/2}(\tilde y - \theta)) p(\tilde y) \mathrm{d}\tilde y}{ \Phi( n_1^{1/2}(\theta - t) )}\right\} \\
& + \Phi\left\{ \frac{n_1^{1/2}(y - t)}{(1-\gamma)^{1/2}}  \right\}.
\end{split}
\end{equation} 
The derivation of this expression can be found in Appendix \ref{APP: pmp} and is the one employed in the numerical implementations considered in the next section.

Since evaluation of $\pi_y(\theta)$ in the $\gamma < 1$ case requires approximating an integral and dividing by small numbers, it can sometimes be numerically unstable. As a computationally lighter alternative, we will also consider the Jeffreys prior of the selective model, given by 
\begin{equation}
\pi_J(\theta) \propto \left\{ n + \frac{\partial^2}{\partial \theta^2} \log \varphi(\theta) \right\}^{1/2}.
\end{equation}
For the one-sided truncation models considered before, this prior can be written as 
\begin{equation} \label{EQ: JP_normal}
\pi_J(\theta) \propto \left[ 1 + \gamma  h_2\{n_1^{1/2}(\theta - t)\} \right]^{1/2}, 
\end{equation}
where $h_2(x) = (\partial^2 / \partial x^2) \log \Phi(x)  =  -x h_1(x) - h_1(x)^2$.

\begin{figure}[H]
\centering
\includegraphics[width = \textwidth]{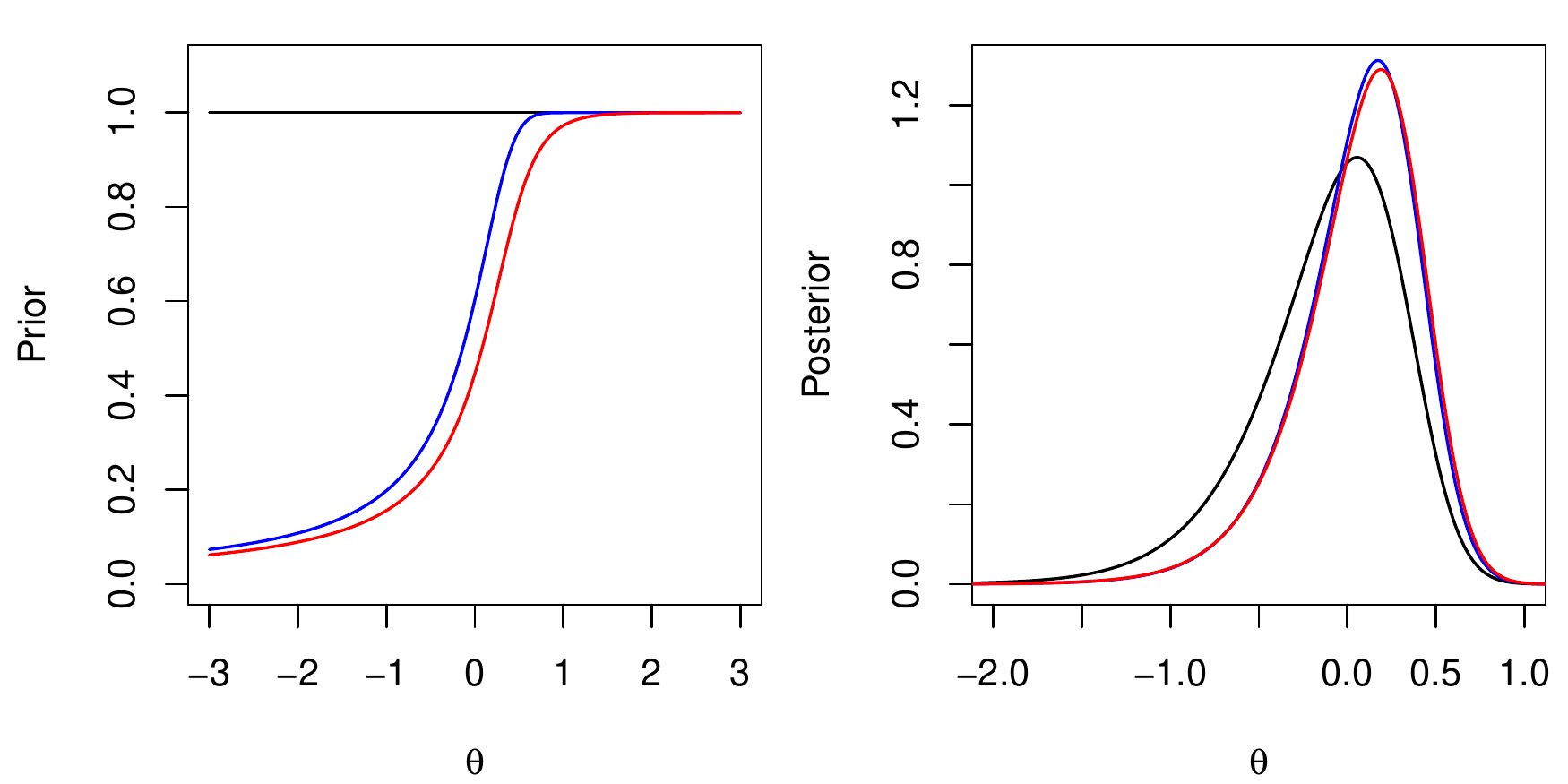}
\caption{Left: uniform prior (black), probability-matching prior for $y = 0.2$ (red), and Jeffreys prior (blue) for the normal model with $(n, \gamma, t) = (20, 1, 0)$. Right: the resulting posteriors for $y = 0.2$ (the red density overlaps the blue one in this panel).}
\label{FIG: noninf1}
\end{figure}

\begin{figure}[H]
\centering
\includegraphics[width = \textwidth]{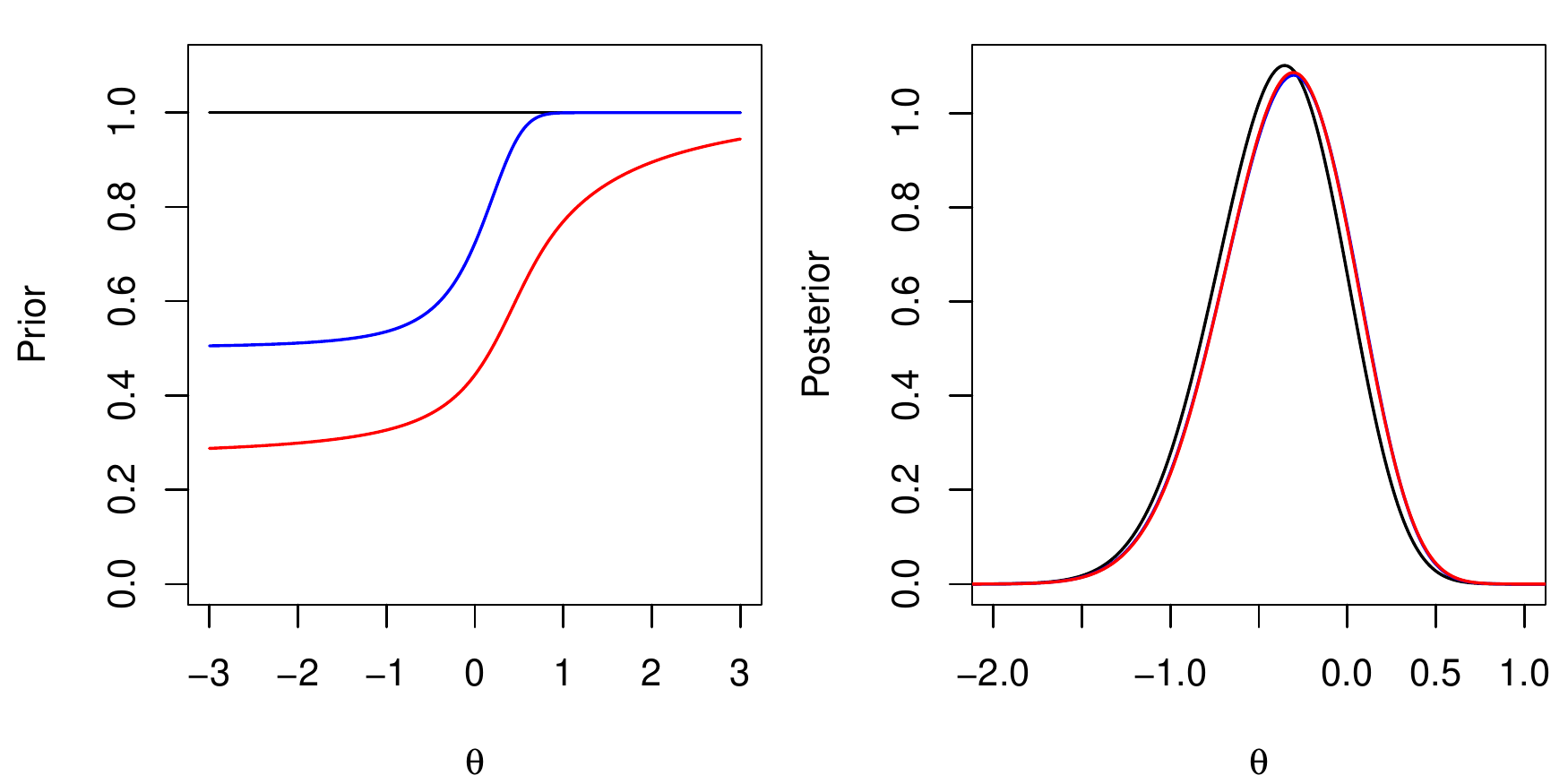}
\caption{Left: uniform prior (black), probability-matching prior for $y = 0$ (red), and Jeffreys prior (blue) for the normal model with $(n, \gamma, t) = (20, 0.75, 0)$. Right: the resulting posteriors for $y = 0$ (the red density overlaps the blue one in this panel).}
\label{FIG: noninf2}
\end{figure}

Note that both priors depend on the sample size, which is unusual in non-selective settings. In Figure \ref{FIG: noninf1} we have plotted both priors and the resulting posterior densities for this model, with $n = 20$, $\gamma = 1$, $t = 0$, and $y = 0.2$, and in Figure \ref{FIG: noninf2} for $n = 20$, $\gamma = 0.75$, $t = 0$, and $y = 0$. The uniform prior and the respective posterior densities are also plotted for comparison. The posterior densities corresponding to the two proposed priors are virtually identical, a behaviour that extends to other choices of parameters and observations, which suggests that Jeffreys prior is well-calibrated from a frequentist viewpoint. Formal theoretical frequentist guarantees are currently unresolved, but we provide empirical evidence of them in Table \ref{TAB: Jeffreys}. The table shows the coverage of the intervals $(-\infty, \Pi^{-1}(\alpha\vert Y)]$ for this model with $n = 20$, $t = 0$, and different combinations of $(\theta, \gamma, \alpha)$, computed via numerical integration. The posterior distributions were constructed using the uniform prior (U) and Jeffreys prior (J). We see that the latter option performs better in most cases, and is considerably superior to the uniform prior when all the data is used for selection. A key characteristic of these priors is that they assign higher prior probabilities to regions of the parameter space with large selection probability $\varphi(\theta)$, thereby correcting the overcompensation of the selection bias described at the start of the section.

\begin{table}[H] 
\centering
\scalebox{0.86}{
\begin{tabular}{ l  *{7}{ c c }  }
\toprule
$\alpha$ & \multicolumn{2}{c }{$0.05$} & \multicolumn{2}{c }{$0.1$} & \multicolumn{2}{c }{$0.25$} & \multicolumn{2}{c }{$0.5$} & \multicolumn{2}{c }{$0.75$} & \multicolumn{2}{c }{$0.9$} & \multicolumn{2}{c }{$0.95$} \\
\toprule
Prior & U & J & U & J & U & J & U & J & U & J & U & J & U & J \\
\toprule
\multicolumn{15}{ c }{$\theta = -0.5$} \\
\hline
$\gamma = .5$ & .045 & .049 & .092 & .099 & .237 & .250 & .487 & .501 & .742 & .751 & .896 & .901 & .948 & .950 \\ 
 $\gamma = .75$ & .035 & .047 & .075 & .096 & .209 & .247 & .458 & .502 & .724 & .754 & .889 & .903 & .944 & .952 \\ 
 $\gamma = 1$ & .020 & .051 & .041 & .098 & .110 & .234 & .250 & .465 & .448 & .714 & .637 & .879 & .738 & .938 \\ 
\hline
\multicolumn{15}{ c }{$\theta = 0$} \\
\hline
$\gamma = .5$ & .043 & .049 & .088 & .099 & .228 & .248 & .473 & .498 & .730 & .750 & .890 & .901 & .945 & .950 \\ 
   $\gamma = .75$ & .038 & .050 & .077 & .099 & .199 & .243 & .427 & .492 & .693 & .748 & .872 & .901 & .935 & .952 \\ 
   $\gamma = 1$ & .033 & .056 & .065 & .108 & .161 & .251 & .331 & .477 & .537 & .710 & .709 & .868 & .794 & .929 \\ 
\hline
\multicolumn{15}{ c }{$\theta = 0.5$} \\
\hline
  $\gamma = .5$ & .048 & .051 & .095 & .101 & .238 & .251 & .479 & .498 & .728 & .747 & .886 & .897 & .941 & .948 \\ 
  $\gamma = .75$ & .048 & .052 & .095 & .103 & .236 & .255 & .471 & .501 & .713 & .744 & .871 & .893 & .930 & .945 \\ 
  $\gamma = 1$ & .048 & .052 & .096 & .105 & .237 & .259 & .469 & .509 & .702 & .749 & .852 & .890 & .908 & .939 \\ 
  \hline
\end{tabular}}
 \caption{\label{TAB: Jeffreys} Estimated coverages of $(-\infty, \Pi^{-1}(\alpha\vert Y)]$ for the normal-location model derived from the uniform (U) and Jeffreys (J) priors.}
\end{table}

\newpage
\subsection{Non-informative priors for exponential families} \label{SEC: PMP_1d}

The derivation of the probability-matching prior in the normal model is only useful to the extent that it can be used to devise priors for other models. If the model is normal we can provide inference directly via the confidence distribution. In this section we use the probability-matching and Jeffreys priors discussed in the previous section to derive non-informative priors for selected one-dimensional parameters of an exponential family, given a random sample from it. This allows us to extend inferential procedures for selective normal location models to other models asymptotically. An attractive feature of exponential families that is useful here is that they always admit a sufficient statistic of the same dimension as the parameter. They are, in fact, the only regular models for which this holds, by the Pitman–Koopman–Darmois Theorem. This enables us to reduce the dimensionality of the selective model to that of the parameter, leading to a substantial mathematical simplification of the problem.

First, consider the one-dimensional case. Suppose that we have a random sample $Y_1, \ldots, Y_n$ from a full exponential family, with density or mass function
\begin{equation}
f(y_i; \theta) = h(y_i)\exp\left\{ \eta(\theta) s(y_i)  - A(\theta) \right\}.
\end{equation}
We assume that the MLE of the non-selective model, $\hat\theta = \argmax_{\theta\in \Theta}  \eta(\theta) \sum_{i = 1}^n s(y_i)  - n A(\theta) $, exists for all realisations of the data. Assume, as before, that the sample is divided into two sets of sizes $n_1$ and $n_2$, $Y_1, \ldots, Y_{n_1}$, and $Y_{n_1 + 1}, \ldots, Y_n$, with respective maximum likelihood estimators $\hat\theta_1$ and $\hat{\theta}_2$, and that selection only uses the first set of samples, so that the selection function may be written as $p(y_1, \ldots, y_{n_1})$. Furthermore, let us denote by $p_{\hat\theta_1}(\hat\theta_1) = \mathbb{E}[p(Y_1, \ldots, Y_{n_1})\mid \hat\theta_1]$ the selection function in terms of $\hat\theta_1$, which is independent of $\theta$, as $\hat\theta_1$ is a sufficient reduction of the first set of samples.

Heuristically, a non-informative prior density for $\theta$ may be derived as follows. Let $i(\theta) = -\mathbb{E}_\theta[(\partial^2/\partial \theta^2)\log f(Y_i;\theta)]$ be the per-sample Fisher information of the non-selective model. Any function $g(\theta)$ satisfying $g'(\theta) = i(\theta)^{1/2}$ is known as a variance-stabilising transformation, as it satisfies 
\begin{equation}
n^{1/2}\{g(\hat{\theta}) - g(\theta)\}\xrightarrow{d} N(0, 1)
\end{equation}
by the Delta Method. In particular, for $i = 1, 2$, we have that $n_i^{1/2}\{g(\hat{\theta}_i) - g(\theta)\}\xrightarrow{d} N(0, 1)$ as $n_i\to \infty$. Thus, for large values of $n_1$ and $n_2$, the original selective model is approximately equivalent to a selective normal-location model involving two independent observations, $U_1\sim N(g(\theta), n_1^{-1})$ and $U_2\sim N(g(\theta), n_2^{-1})$, and selection function $p(u_1, u_2) = p_{\hat\theta_1}\{g^{-1}(u_1)\}$. This argument suggests choosing as a non-informative prior for the parameter $\nu = g(\theta)$ a prior that is reliable in the limiting normal model. Denoting the chosen prior by $\pi_\nu(\nu)$, the resulting prior in the original parametrisation is given by $\pi_\theta(\theta) \propto i(\theta)^{1/2} \pi_\nu\{g(\theta)\}$. In view of the previous discussion and the simulation results, we will take $\pi_\nu(\nu)$ to be either the probability-matching prior or the Jeffreys prior. 
Note that, in the absence of selection, both choices give $\pi_\nu(\nu)\propto 1$, so $\pi_\theta(\theta)$ specialises to the Jeffreys prior of the non-selective model, $\pi_\theta(\theta)\propto i(\theta)^{1/2}$, which is probability-matching to order $O(n^{-1})$.

\begin{figure}[ht]
\includegraphics[width=\textwidth]{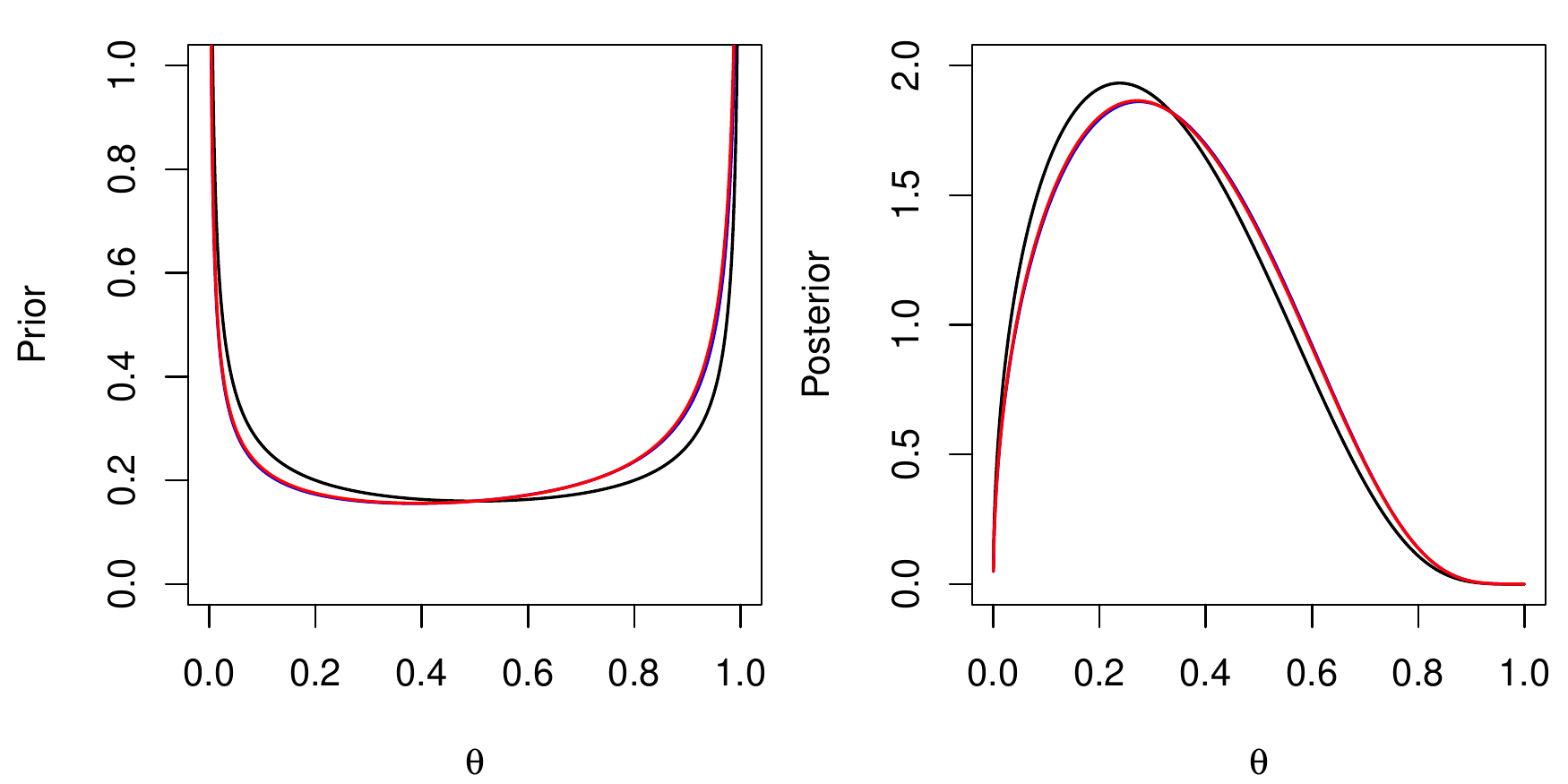}
\caption{Left: non-selective Jeffreys prior (black), probability-matching prior for $y_1 = 4$ and $y_2 = 1$ (red), and Jeffreys prior (blue) for the binomial model. Right: the resulting posteriors (the red lines overlap the blue ones).}
\label{FIG: binomial}
\end{figure}

As an example, consider a model involving two binomial observations $Y_1 \sim \mathrm{Bin}(n_1, \theta)$ and $Y_2 \sim \mathrm{Bin}(n_2, \theta)$ and selection event $n_1^{-1}Y_1 \geq 0.5$, where $\theta\in [0, 1]$. The variance-stabilising transformation of this model is $g(\theta) = \sin^{-1}(\theta^{1/2})$. Therefore, the option based on Jeffreys prior produces 
\begin{equation}
\pi_J(\theta) \propto \theta^{-1/2}(1-\theta)^{-1/2} \left[ 1 + \gamma  h_2 \left\{n_1^{1/2}\left(\sin^{-1}(\theta^{1/2}) - \frac{\pi}{4}\right)\right\} \right]^{1/2},
\end{equation}
and similarly for the option based on the exact probability-matching prior, $\pi_y(\theta)$, which we do not reproduce here. The left panel of Figure \ref{FIG: binomial} shows the standard, non-selective Jeffreys prior $\pi(\theta)\propto \theta^{-1/2}(1-\theta)^{-1/2}$, and the two non-informative selective priors for $n_1 = 8$, $n_2 = 2$, $y_1 = 4$, and $y_2 = 1$. The right panel shows the corresponding posterior densities. As in the normal case, the two selective priors are almost identical and favour parameter values with large selection probability.

In the following examples we illustrate the repeated-sampling performance of the proposed priors, showing that they do in fact lead to well calibrated posterior inference in a frequentist sense. 

\begin{figure}[H]
\centering
\includegraphics[width = \textwidth, height = 1.2\textwidth]{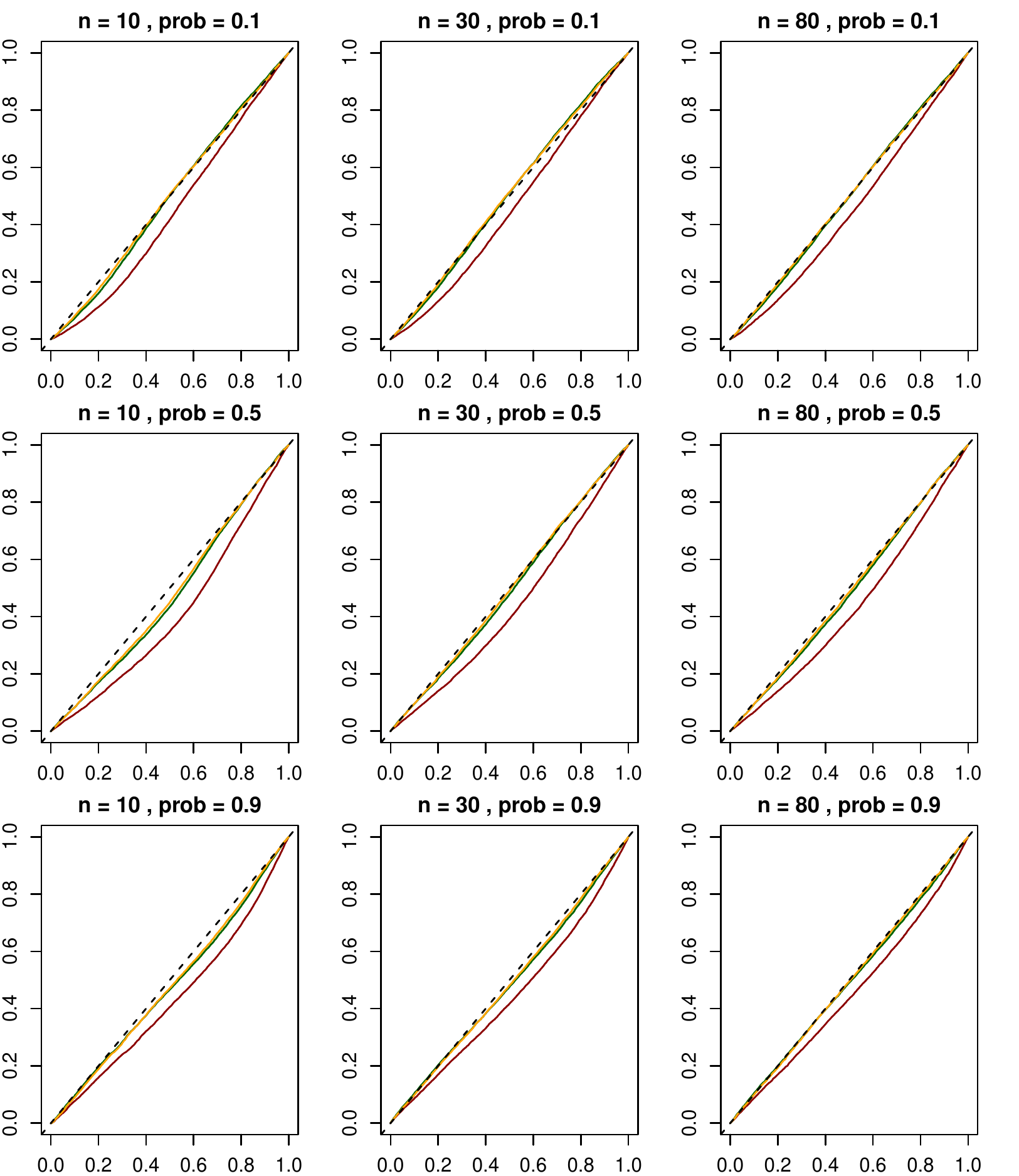}
\caption{Example \ref{EX: exponential}. Coverage of the interval $(-\infty, \Pi^{-1}(\alpha\vert Y_1, \ldots, Y_n)]$ as a function of $\alpha$ for the non-selective Jeffreys prior (red), the probability-matching prior (orange), and the selective Jeffreys prior (green). The orange lines partially overlap the green ones.}
\label{FIG: exponential}
\end{figure}

\begin{example} \label{EX: exponential}
\textit{Exponential}.
Let $Y_1, \ldots, Y_n$ be a random sample from an exponential distribution with rate parameter $\theta > 0$. For this model the variance-stabilising transformation is $g(\theta) = \log(\theta)$. We consider the selection event $\hat{\theta}_1 > 1$, where $\hat{\theta}_1$ is the maximum likelihood estimator of $\theta$ based on a subsample of size $n_1 = 0.8\times n$. For the simulation we consider the values of $n = 10, 30, 80$, and for each sample size we consider three true parameter values $\theta_0$ defined to satisfy $\varphi_n(\theta_0) = 0.1, 0.5, 0.9$, corresponding to situations with small, moderate and large selection probabilities, respectively. For each pair $(n, \theta_0)$, we plot the coverage of the interval $(-\infty, \Pi^{-1}(\alpha\vert Y_1, \ldots, Y_n)]$ as a function of $\alpha$ for the non-selective Jeffreys prior $\pi(\theta)\propto \theta^{-1}$ and for the two non-informative priors proposed in this work. The coverages were approximated via $10^4$ simulations from the conditional model $(Y_1, \ldots, Y_n)^T\vert \hat{\theta}_1 > 1$. The results can be found in Figure \ref{FIG: exponential}. The performances of the selective Jeffreys and probability-matching priors are practically identical and significantly better than that of the non-selective prior.
\end{example}

\begin{figure}[H]
\centering
\includegraphics[width = \textwidth, height = 1.2\textwidth]{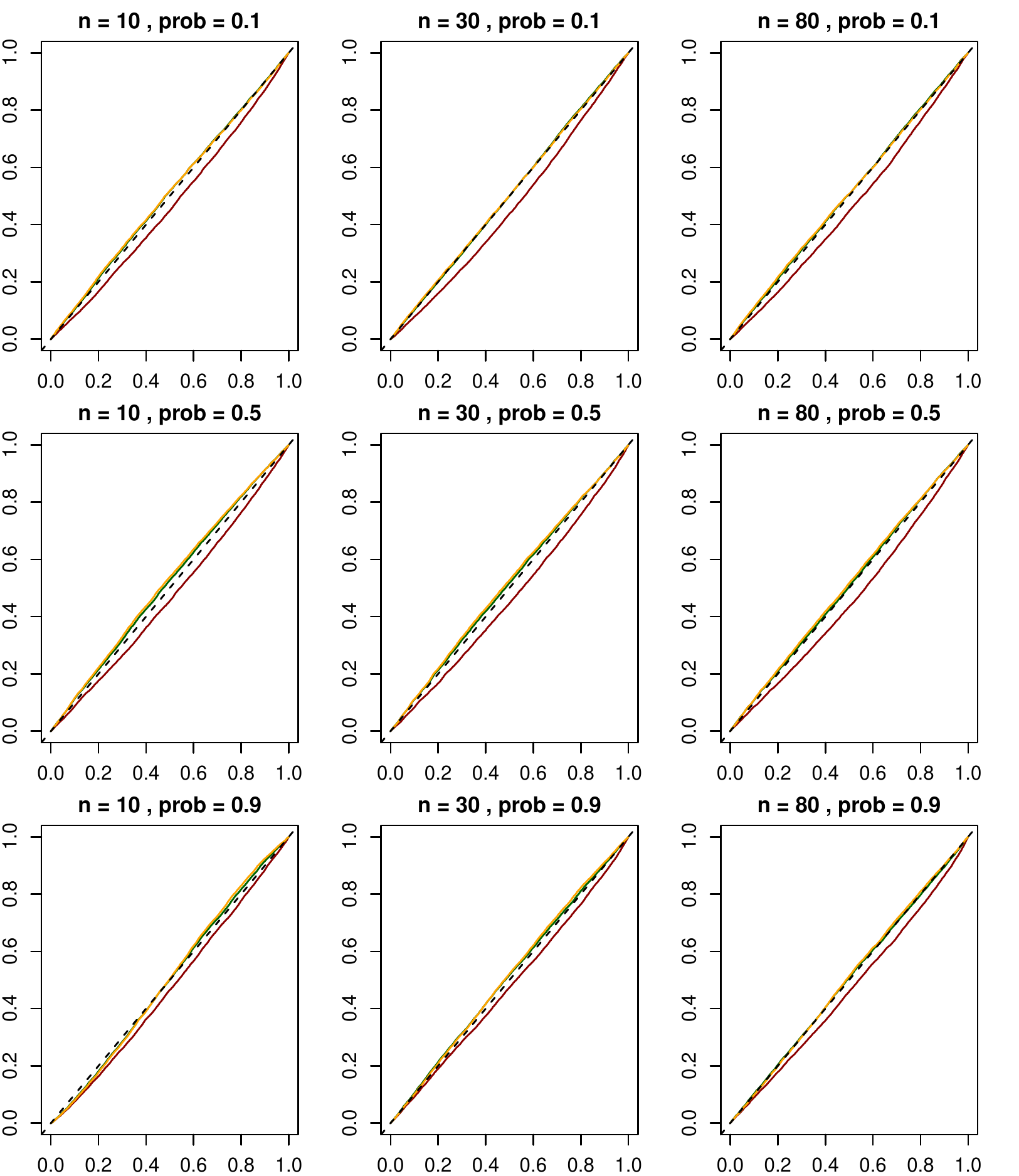}
\caption{Example \ref{EX: invgauss}. Coverage of the interval $(-\infty, \Pi^{-1}(\alpha\vert Y_1, \ldots, Y_n)]$ as a function of $\alpha$ for the non-selective Jeffreys prior (red), the probability-matching prior (orange), and the selective Jeffreys prior (green). The orange lines partially overlap the green ones.}
\label{FIG: invgauss}
\end{figure}

\begin{example} \label{EX: invgauss}
\textit{Inverse Gaussian}.
Let $Y_1, \ldots, Y_n$ be a random sample from an inverse Gaussian distribution with unknown mean $\theta > 0$ and known shape parameter $\lambda = 1$, with density
\begin{equation}
f(y_i; \theta) = \frac{1}{(2\pi y_i^3)^{1/2}}\exp\left\{ -\frac{(y_i - \theta)^2}{2\theta^2 y_i} \right\}.    
\end{equation}
The variance-stabilising transformation for this model can be easily seen to be $g(\theta) = -2\theta^{-1/2}$, and the non-selective Jeffreys prior $\pi(\theta)\propto \theta^{-3/2}$. We consider the same settings as before. The selection event is $\hat{\theta}_1 > 1$, where $\hat{\theta}_1$ is the maximum likelihood estimator of $\theta$ based on a subsample of size $n_1 = 0.8\times n$ for $n = 10, 30, 80$, and we define the true values of the parameter in the same way as before. The results, shown in Figure \ref{FIG: invgauss}, are similar to those of the exponential model. 
\end{example}

We now consider situations where the full parameter $\theta$ is multidimensional and the selected parameter, $\psi$, is one-dimensional. Again, we assume that the non-selective model is a full exponential family, with density or mass function
\begin{equation}
f(y_i; \theta) = h(y_i)\exp\left\{ \eta(\theta)^T s(y_i) - A(\theta) \right\}, \quad \theta\in\Theta\subseteq \mathbb{R}^p,
\end{equation}
and, as before, we assume that selection is based on a random sample of size $n_1$, $Y_1,\ldots, Y_{n_1} $, and that there is a second set of samples of size $n_2$, $Y_{n_1 + 1}, \ldots, Y_n$, available for inference. 

Suppose that we can write $\theta = (\psi, \chi)$, and let us reparametrise the model to $(\psi, \lambda)$, where $\lambda \equiv \lambda(\psi, \chi)$ is orthogonal to $\psi$ in the sense of \cite{orthogonal}. This means that the Fisher information matrix in the latter parametrisation has zeroes in all the non-diagonal entries of the first row and column, or equivalently, that the maximum likelihood estimators $\hat\psi$ and $\hat\lambda$ are asymptotically independent, in the absence of selection. Such a parametrisation always exists when $\psi$ is one-dimensional, and admits a very simple expression in terms of the cumulant function $A(\theta)$ if the model its parametrised in its natural form, so that $\eta(\theta) = \theta$. 

Let us denote by $i_{\psi\psi}(\psi, \lambda)$ and $i_{\lambda\lambda}(\psi, \lambda)$ the components of the Fisher information matrices corresponding to $\psi$ and $\lambda$ respectively. Also, let $\hat\psi_i$ and $\hat\lambda_i$, $i = 1, 2$, be the maximum likelihood estimators based on the first and second set of samples, assumed to exist, and let $p_{\hat\psi_1, \hat\lambda_1}(\hat\psi_1, \hat\lambda_1)$ be the selection function in terms of $(\hat\psi_1, \hat\lambda_1)$. By orthogonality, the asymptotic distribution of $\hat\psi_1$ given $\hat\lambda_1$ and the selection event is selective normal with mean $\psi$, variance $\{n_1 i_{\psi\psi}(\psi, \lambda)\}^{-1}$, and selection function $p_{\hat\psi_1}(\hat\psi_1) = p_{\hat\psi_1, \hat\lambda_1}(\hat\psi_1, \hat\lambda_1)$, where $\hat\lambda_1$ is fixed at its observed value. Also, the asymptotic distribution of $\hat\psi_2$ given $\hat\lambda_2$ is $N(\psi, \{n_2 i_{\psi\psi}(\psi, \lambda)\}^{-1})$. Therefore, for a fixed value of $\lambda$, $(\hat\psi_1, \hat\psi_2)\vert (\hat\lambda_1, \hat\lambda_2)$ follows asymptotically a selective normal distribution, with variance-stabilising transformation $g(\psi; \lambda)$ satisfying $g'(\psi; \lambda) = i_{\psi\psi}(\psi, \lambda)^{1/2}$. Since $\hat\lambda_1$ is only mildly informative about $\psi$, the same heuristic argument as before would lead to adoption of the conditional prior $\pi(\psi\vert \lambda) \propto  i_{\psi\psi}(\psi, \lambda)^{1/2}  \pi_\nu\{g(\psi; \lambda)\}$, where $\pi_\nu(\nu)$ is a suitable non-informative prior for the selective normal-location model with location parameter $\nu = g(\psi; \lambda)$. Since this choice does not constrain the marginal prior of $\lambda$, a generic joint prior for $(\psi, \lambda)$ would be of the form 
\begin{equation}
\pi(\psi, \lambda)\propto c(\lambda) i_{\psi\psi}(\psi, \lambda)^{1/2}  \pi_\nu\{g(\psi; \lambda)\}, 
\end{equation}
where $c(\lambda)$ is an arbitrary prior density for $\lambda$.

\newpage
\subsubsection{Inference for a selected normal mean with an unknown variance} \label{SEC: normal unknown}

A simple selection problem involving nuisance parameters occurs when inference is required for a selected normal mean and the variance is unknown. Let $Y_1, \dots, Y_n\sim N(\mu, \sigma^2)$, where both $\mu$ and $\sigma^2$ are unknown, but only $\mu$ is of direct interest to us. Suppose that the sample is divided into two sets of sizes $n_1$ and $n_2$, with respective maximum likelihood estimators $\hat{\theta}_1 = (\bar{Y}_1, V_1)$ and $\hat{\theta}_2 =(\bar{Y}_2, V_2)$ of $\theta = (\mu, \sigma^2)$. Suppose that, in order to determine whether $\mu$ is of interest, we conduct the $t$-test $ V_1^{-1/2}\bar{Y}_1 > n_1^{-1/2}t$ for some pre-specified value of $t$. The selective density of the data is
\begin{equation}\label{EQ: sel_den_unkownvar}
f_S(y_1, \ldots, y_n; \mu, \sigma^2) = \frac{f(y_1, \ldots, y_n; \mu, \sigma^2)\mathbf{1}(v_1^{-1/2}\bar{y}_1 > n_1^{-1/2}t)}{\mathbb{P}_{\mu, \sigma^2}(V_1^{-1/2}\bar{Y}_1 > n_1^{-1/2}t)}.
\end{equation}
Note that in the selective model the marginal distribution of $V_1$, with marginal density
\begin{equation} \label{EQ: sel_den_v1}
f_S(v_1; \mu, \sigma^2) = \frac{f(v_1; \sigma^2)\mathbb{P}_{\mu, \sigma^2}(v_1^{-1/2}\bar{Y}_1 > n_1^{-1/2}t\mid v_1)}{\mathbb{P}_{\mu, \sigma^2}(V_1^{-1/2}\bar{Y}_1 > n_1^{-1/2}t)},
\end{equation}
depends on $\mu$, even though its non-selective density is free of it. However, in this case the non-selective distribution of $\bar{Y}_1\vert V_1$ is not free of $\sigma^2$, so there is no justification for conditioning on $V_1$.

\begin{figure}[H]
\centering
\includegraphics[scale = .6]{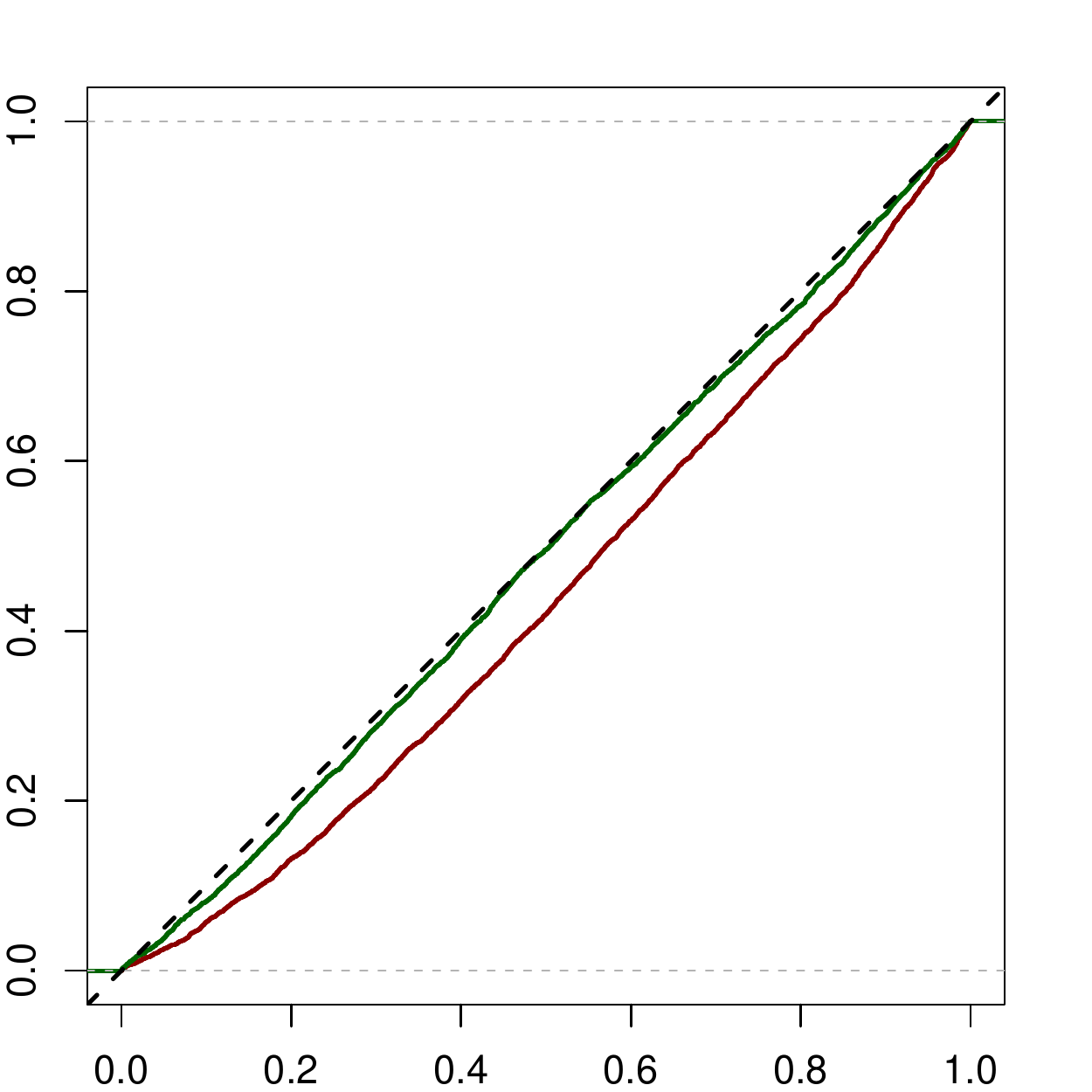}
\caption{Estimated empirical CDFs of $\Pi(\mu_0\vert Y_1, \ldots, Y_n)$ for the prior $\pi(\mu, \sigma^2)\propto \sigma^{-1}$ (red) and for selective Jeffreys prior (green).}
\label{FIG: 4}
\end{figure}

As expected by comparison with the univariate case, standard non-informative priors such as $\pi(\mu, \sigma^2)\propto \sigma^{-1}$ produce marginal posteriors for $\mu$ that overstate, on average, smaller values of $\mu$. In this case, the proposed non-informative priors can be found easily, as the model is already parametrised orthogonally. For example, the option based on the Jeffreys prior takes the form
\begin{equation}
\pi_J(\mu, \sigma^2; v_1) \propto c(\sigma^2) \left\{ 1 + \frac{n_1}{n} h_2\left( \frac{n_1^{1/2}\mu}{\sigma} - \frac{t v_1^{1/2}}{\sigma}\right) \right\}^{1/2},
\end{equation}  
where $v_1$ is the observed value of $V_1$, and the option based on the probability-matching prior for $\gamma = 1$ is
\begin{equation}
\pi_y(\mu, \sigma^2; v_1) \propto c(\sigma^2) \left\{ 1 - \frac{h_1 ( \sigma^{-1} n^{1/2}\mu - \sigma^{-1}t v^{1/2})}{h_1 ( \sigma^{-1} n^{1/2}\mu - \sigma^{-1} n^{1/2}\bar y)}  \right\}.
\end{equation} 
The corresponding expression for the case $\gamma < 1$ is rather complicated and therefore omitted, but can be obtained easily from \eqref{EQ: PMP_normal_long}. These priors assign lower probabilities to small values of $\mu$ relative to $\sigma$. Figure \ref{FIG: 4} illustrates the performance of $\pi_J(\mu, \sigma^2; v_1) $ for $c(\sigma^2) = \sigma^{-1}$ under repeated sampling from the selective model with $n_1 = 50$, $n_2 = 10$, $t = 2$, and true parameter $(\mu_0, \sigma^2_0) = (0, 1)$. The prior $\pi_y(\mu, \sigma^2; v_1) $ gave numerical issues for small values of $\sigma$ and was not considered. The plot shows the empirical distribution functions of $\Pi(\mu_0\vert Y_1, \ldots, Y_n)$ for the prior $\pi(\mu, \sigma^2)\propto \sigma^{-1}$ and for the proposed prior, as estimated from $5\times 10^3$ repetitions. The selection probability was computed by numerical evaluation of the integral
\begin{equation}
\mathbb{P}_{\mu, \sigma^2}(V_1^{-1/2}\bar{Y}_1 > n_1^{-1/2}t) = \int_0^\infty f(v_1;\sigma^2) \Phi\left\{ \frac{n_1^{1/2}}{\sigma}\left( \mu - \frac{t v_1^{1/2}}{n_1^{1/2}} \right) \right\} \mathrm{d}v_1,
\end{equation} 
and the marginal posterior distribution of $\mu$ was approximated with a Metropolis-Hastings algorithm with $5\times 10^3$ steps. The results show that the selection-adjusted prior produces posterior inference with a more reliable frequentist calibration than the unadjusted one.

\subsubsection{Inference for the winner} \label{SEC: play winner}

Consider the standard problem of providing inference for the mean of a largest normal sample. Suppose we have $n_1$ observations from $m$ different means, with sample means  $Y_i\sim N(\theta_i, n_1^{-1})$, $i = 1, \ldots, m$, and assume without loss of generality that $Y_1 = \max\{Y_i\colon i = 1, \ldots, m\}$, so that $\theta_1$ is the parameter of interest. Assume also that after selection we obtain another sample of size $n_2$ from $\theta_1$, with mean $\tilde{Y}_1 \sim N(\theta_1, n_2^{-1})$. There are two natural inferential models that are compatible with this selection rule. One has likelihood
\begin{equation}
L_1(\theta) \propto \frac{\phi\{n_2^{1/2}(\theta_1 - \tilde y_1)\} \prod_{i = 1}^m \phi\{n_1^{1/2}(\theta_i - y_i)\}}{\mathbb{E}_\theta[\Phi\{n_1^{1/2} (\theta_1 - T) \}]}
\end{equation}
where $n = n_1 + n_2$, $\theta = (\theta_1, \ldots, \theta_m)^T$, and $T = \max\{Y_i\colon i = 2, \ldots, m\}$, and the other has likelihood
\begin{equation}
L_2(\theta) \propto \frac{\phi\{n_2^{1/2}(\theta_1 - \tilde y_1)\} \prod_{i = 1}^m \phi\{n_1^{1/2}(\theta_i - y_i)\}}{\Phi\{n_1^{1/2} (\theta_1 - t) \}}.
\end{equation}
The latter is derived from the reasonable procedure of conditioning on the data $Y_2, \ldots, Y_m$ from the non-selected means: these only depend mildly
on the parameter of interest, under selection, and the resulting inference is simpler, reducing to a one-dimensional normal problem, as we saw before.

For this model, the prior proposed in the discussion above takes the form
\begin{equation} \label{EQ: winner_prior}
\pi(\theta) \propto c(\theta_2, \ldots, \theta_m) \left[  1 + \frac{n_1}{n} h_2\left\{  n_1^{1/2} (\theta_1 - t)  \right\} \right]^{1/2}
\end{equation}
for the case of utilising the Jeffreys prior, and can similarly be constructed based on the probability-matching prior. Note that in this case the original parametrisation is also orthogonal.

Since the distribution of $Y_2, \ldots, Y_m$ depends on $\theta_1$ given selection, conditioning on them, even if advisable on computational grounds, will generally reduce the inferential power. In Table \ref{TAB: winner} we investigate this by comparing the length of credible intervals obtained by attaching prior \eqref{EQ: winner_prior} with $c(\theta_2, \ldots, \theta_m) = 1$ to the likelihoods $L_1(\theta)$ and $L_2(\theta)$. We set $n_1 = n_2= 5$, true parameter $\theta = (0, \ldots, 0)^T$, and $m = 2, 5, 10, 20$. The credible intervals were computed by taking the lower and upper $0.05$ quantiles of the marginal posterior distribution of $\theta_1$. To compute each of the intervals we run a Metropolis-Hastings algorithm with $10^4$ steps. The simulation size was $5\times 10^3$, and the maximum standard error of the figures in the table was $0.005$. We can see there is some loss of power when we condition on the data from the non-selected means, as indicated by longer intervals, and that this loss worsens slightly as the number of means increases. We note, however, that the simulation for $L_1(\theta)$ was much more computationally challenging than for $L_2(\theta)$ for large values of $m$. The results for $L_1(\theta)$ were found to be quite sensitive to the choice of integration boundaries in the approximation of the selection probability. Similarly, use of the prior based on probability-matching is prone to numerical instability, as it involves a high-dimensional numerical integration. Again, use of the Jeffreys prior is strongly supported in terms of frequentist calibration.

\begin{table}[ht]
\centering
\begin{tabular}{l c c c c}
 & \multicolumn{2}{c}{$L_1(\theta)$} & \multicolumn{2}{c}{$L_2(\theta)$} \\
 \toprule
  & Cov. & Length & Cov. & Length \\
  \toprule
  	$m = 2$ & 0.913 & 1.149 & 0.902 & 1.205 \\
  	$m = 5$ & 0.899 & 1.193  & 0.902 & 1.292  \\
  	$m = 10$ & 0.899  & 1.183 & 0.906 & 1.329 \\
  	$m = 20$ & 0.887 & 1.148  & 0.903 & 1.352  \\
  	\bottomrule
\end{tabular}
\caption{Coverage and average length of 0.9-credible intervals for the mean giving the largest sample, derived from likelihoods $L_1(\theta)$ and $L_2(\theta)$ and Jeffreys prior.}
\label{TAB: winner}
\end{table}

\section{Discussion}

We have analysed the selective inference problem from a Bayesian perspective, arguing in favour of a selection adjustment of the inferences, contrary to the classical viewpoint. Furthermore, we have proposed two classes of non-informative prior densities for selection models which are approximately probability matching. These priors provide some level of regularisation in low-probability regions and behave better than alternatives proposed in the literature which are independent of selection. To derive these priors we first considered the case of normal-location models, and then extended the analysis to full exponential families by appealing to asymptotic considerations. One of the proposed priors provides exact matching in the normal case by construction. The other is the Jeffreys prior, constructed from the selective likelihood, which provides similar results and is more computationally stable. These priors were derived heuristically and their behaviour was demonstrated empirically. Formal analyses are challenging due to the non-standard asymptotic behaviour of selective models. Theory on probability-matching priors typically relies on second-order normal approximations to the posterior which are not applicable in these settings. Further theoretical study is warranted.

\bibliographystyle{agsm}
\bibliography{references}

\vspace{2cm}
{\LARGE\textbf{Appendix}}

\appendix 

\section{Proof of Proposition \ref{PROP: uniform prior}}\label{APP: prop}

Without loss of generality we may assume that $\sigma^2 = 1$ and $t = 0$. Let $h_i(x)$ denote the $i$-th derivative of $\log \Phi(x)$ and, for $y>0$, define
\begin{equation}
H(\theta; y) = \mathbb{P}_\theta(Y\geq y\vert Y> 0) = \frac{\Phi(\theta - y)}{\Phi(\theta)},     
\end{equation}
which is a distribution function as a function of $\theta$, and
\begin{equation}
g(\theta; y) = -\frac{\frac{\partial}{\partial \theta} H(\theta; y)}{\frac{\partial}{\partial y} H(\theta; y)} = 1 - \frac{h_1(\theta)}{h_1(\theta - y)} > 0.
\end{equation}
This function is strictly increasing. Indeed,
\begin{equation}
g'(\theta; y) =\frac{h_1(\theta)h_2(\theta - y)}{h_1(\theta - y)^2} - \frac{h_2(\theta)}{h_1(\theta - y)},
\end{equation}
which is positive if 
\begin{equation}
\frac{h_2(\theta - y)}{h_1(\theta - y)} > \frac{h_2(\theta)}{h_1(\theta)},
\end{equation}
that is, if $h_2(x)/h_1(x)$ is strictly decreasing. Since  $h_2(x) = - xh_1(x) - h_1(x)^2$, $h_2(x)/h_1(x) = - x - h_1(x)$, so $(\partial/\partial x)\{h_2(x)/h_1(x)\} = - 1 - h_2(x) < 0$ for all $x \in \mathbb{R}$. The latter claim follows from the fact that $1 + h_2(x)$ is the variance of a $N(x, 1)$ distribution truncated to $(0, \infty)$. Since $g(\theta; y)$ is strictly increasing, if used as a formal prior density for $\theta$, it satisfies
\begin{equation}
\Pi(\theta\vert y) > \frac{\int_{-\infty}^{\theta_0} g(\theta; y) \Phi(\theta)^{-1}\phi(\theta - y)\mathrm{d}\theta}{\int_{-\infty}^{\infty} g(\theta; y) \Phi(\theta)^{-1}\phi(\theta - y)\mathrm{d}\theta}.
\end{equation}
But, by definition, $g(\theta; y)$ is such that
\begin{equation}
g(\theta; y) \frac{\phi(\theta - y)}{\Phi(\theta)} = \frac{\partial}{\partial \theta} H(\theta; y),
\end{equation}
so
\begin{equation}
\mathbb{P}_{\theta_0}\{\Pi(\theta_0\vert Y) < \alpha \vert Y>0\} < \mathbb{P}_{\theta_0}\{H(\theta_0, Y)<\alpha\vert Y>0\} = \alpha,
\end{equation}
since $H(\theta_0; Y)$ is uniformly distributed given selection.

\section{Derivation of the probability-matching prior with data splitting} \label{SEC: PMP} \label{APP: pmp}
Recall that the selection function for $\gamma < 1$ can be written as $\mathbb{P}(y + W > t)$, where $W\sim N(0, n^{-1}\{\gamma^{-1} - 1\})$, so
\begin{align}
p(y) =  \Phi\left\{ \left( \frac{n}{\gamma^{-1}-1} \right)^{1/2}  (y-t)\right\}.
\end{align}
We also have that
\begin{equation}
\varphi(\theta) = \Phi\{n_1^{1/2}(\theta - t) \}.
\end{equation}
The definition gives
\begin{equation} \label{EQ: PMP_proof}
\begin{split}
\pi_y(\theta) \propto \frac{1}{\phi(n^{1/2}(y - \theta))}\bigg\{ \int_y^\infty n(\tilde y-\theta)\phi(n^{1/2}(\tilde y - \theta)) p(\tilde y) \mathrm{d}\tilde y & \\ 
 - \frac{\varphi'(\theta)}{\varphi(\theta)} \int_y^\infty \phi(n^{1/2}(\tilde y - \theta)) p(\tilde y) \mathrm{d}\tilde y & \bigg\} .
\end{split}
\end{equation}
Using that
\begin{equation}
\int x \phi(x) \Phi(a + bx) \mathrm{d}x = \frac{b}{d}\phi\left( \frac{a}{d}\right)\Phi\left( xt + \frac{ab}{d} \right) - \phi(x) \Phi(a + bx)
\end{equation}
for any constants $a$ and $b$ \citep{Owen_integrals}, where $d = (1 + b^2)^{1/2}$, the first integral can be expressed as
\begin{equation}
\begin{split}
\gamma^{1/2} \phi(n_1^{1/2}(\theta - t))\left\{ 1 - \Phi\left\{ \frac{n^{1/2}}{(1 - \gamma)^{1/2}}(y - \theta + \gamma(\theta - t)) \right\}\right\} & \\
+ \phi(n^{1/2}(y - \theta))  \Phi\left\{ \frac{n_1^{1/2}}{(1 - \gamma)^{1/2}}(y - t) \right\}& .
\end{split}
\end{equation}
Substituting this expression in \eqref{EQ: PMP_proof} and reorganising the terms gives the claimed expression.

\end{document}